\documentclass[]{bolmat5}
\usepackage{amsfonts}
\usepackage{amssymb}
\usepackage{amsmath}
\usepackage{tikz-cd}
\usepackage{float}
\usepackage{tikz}
\usetikzlibrary{trees}
\usetikzlibrary{decorations.pathmorphing}
\usetikzlibrary{decorations.pathreplacing}
\usetikzlibrary{decorations.markings}
\usetikzlibrary{decorations.shapes}
\usepackage{verbatim}
\usepackage{amsmath}
\usepackage{amsthm}
\usepackage{dsfont}
\usepackage[T1]{fontenc}
\usepackage{babel}[spanish]
\usepackage{makeidx}
\usepackage{enumerate}
\usepackage{calc}
\usepackage{latexsym}
\usepackage{amssymb}
\usepackage{amsfonts}
\usepackage[mathscr]{eucal}
\usepackage{amscd}
\usepackage{array}
\usepackage{amsrefs}
\usepackage{tikz}
\usepackage{listings}
\usepackage{chngcntr}
\usepackage{setspace}

\counterwithin*{equation}{section}
\counterwithin*{equation}{subsection}

\newcommand{\xitem}{%
  \par\hangindent3em\hangafter1
  \noindent\makebox[3em][l]{$\triangleright$}%
  \ignorespaces}

\labeldocument[firstpage = 1, volume = 0, number = 0, month = 00, year = 1900, day = 00, monthreceived = 0, yearreceived = 1900, monthaccepted = 0, yearaccepted = 1900]

\begin{document}
	
\title[maintitle = {Comparing Lagrange and Mixed finite element methods using MFEM library},
	othertitle = {Comparando los m\'etodos de elementos finitos de Lagrange y mixto utilizando la librer\'ia MFEM},
	shorttitle = {}
]
	
\begin{authors}[] 
\author[firstname = {Felipe},
	surname = {Cruz},
	institutionnumber = {1},
	email = {fcruzv@unal.edu.co},
]
\end{authors}
	
\begin{affiliations}
\affiliation[
	department = {Departamento de Matem\'aticas},
	institution = {Universidad Nacional de Colombia},
	city = {Bogot\'a D.C.},
	country = {Colombia}
]
\end{affiliations}

\begin{mainabstract} In this paper, we develop two finite element formulations for the Laplace problem and find the way in which they are equivalent. Then we compare the solutions obtained by both formulations, by changing the order of the shape functions and the refinement level of the mesh (star with rhomboidal elements). And, we will give an overview of MFEM library from the LLNL (Lawrence Livermore National Laboratory), as it is the library used to obtain the solutions.
		
\end{mainabstract}
	
\begin{otherabstract} En este art\'iculo, desarrollamos dos formulaciones de elementos finitos, la de Lagrange y la mixta, y encontramos la manera en que son equivalentes. Luego, comparamos las soluciones obtenidas mediante ambas formulaciones al cambiar el grado de las "shape functions" y el nivel de refinamiento de la malla (una estrella con elementos romboidales). Y, daremos una revisi\'on general de la librer\'ia MFEM, ya que es la librer\'ia utilizada para obtener las soluciones.
		
		
\end{otherabstract}


\textit{\texttt{Note:} This work was done during the second period of 2020 in the course "Beyond Research" from the National University of Colombia. It was supervised by Juan Galvis and Boyan Lazarov.}
\section{Theoretical framework}
In this section we are going to study the theoretic background of the project. First, we are going to review the two finite element methods used (with the problem they solve) and then, give some information about the library. In the finite element parts we'll develop a problem and define the finite element spaces used; all this in two dimensions. And, for the library part, we'll give an overview of its characteristics and the general structure of the code.

\subsection{Lagrange finite elements}
For this method, we consider the following problem \cite{Lagrange}:
\begin{equation}
\begin{split}
        -\Delta &p=f\text{ in } \Omega\\
        &p=0\text{ in }\Gamma
\end{split}
\end{equation}
where $\Omega\subseteq\mathbb{R}^2$ is an open-bounded domain with boundary $\Gamma$, $f$ is a given function and $\Delta p=\frac{\partial^2p}{\partial x^2}+\frac{\partial^2p}{\partial y^2}$.
Consider the space $V$:
\begin{equation*}
    V=\{v:v \text{ continuous on } \Omega, \frac{\partial v}{\partial x},\frac{\partial v}{\partial y} \text{ piecewise continuous on }\Omega\text{ and } v=0 \text{ on } \Gamma\}
\end{equation*}
Now, we can multiply in the first equation of (1) by some $v\in V$ ($v$ is called \textit{test function}) and integrate over $\Omega$:
\begin{equation}
    -\int_\Omega \Delta p\ v=\int_\Omega f\ v
\end{equation}
Applying divergence theorem, the following Green's formula can be deduced \cite{Lagrange}:
\begin{equation}
    -\int_\Omega \Delta p\ v=\int_\Omega\nabla v\cdot\nabla p-\int_\Gamma v\ \nabla p\cdot\eta
\end{equation}
where $\eta$ is the outward unit normal to $\Gamma$.\\
\newline
Since $v=0$ on $\Gamma$, the third integral equals $0$.\\
\newline
\textit{Remark: The boundary integral does not depend on $p$'s value on $\Gamma$ but rather on it's derivative in $\Gamma$. And, this is what's called an \textbf{essential} boundary condition.}\\
\newline
Then, replacing (3) on (2), we get:
\begin{equation}
    \int_\Omega\nabla v\cdot\nabla p=\int_\Omega f\ v
\end{equation}
\textbf{\textit{Note:}}\cite{Lagrange} If $p\in V$ satisfies (4) for all $v\in V$ and is sufficiently regular, then $p$ also satisfies (1), ie, it's a solution for our problem.\\
\newline
In order to set the problem for a computer to solve it, we are going to discretize it and encode it into a linear system.\\
\newline
First, consider a \textit{triangulation} $T_h$ of the domain $\Omega$. This is, $T_h=\{K_1,\dots,K_m\}$ a set of non-overlapping triangles such that $\Omega=K_1\cup\dots\cup K_m$ and no vertex ($N_i$) of one triangle lies on the edge of another triangle:
\newpage
\begin{figure}[h!]
    \centering
    \includegraphics[scale=0.3]{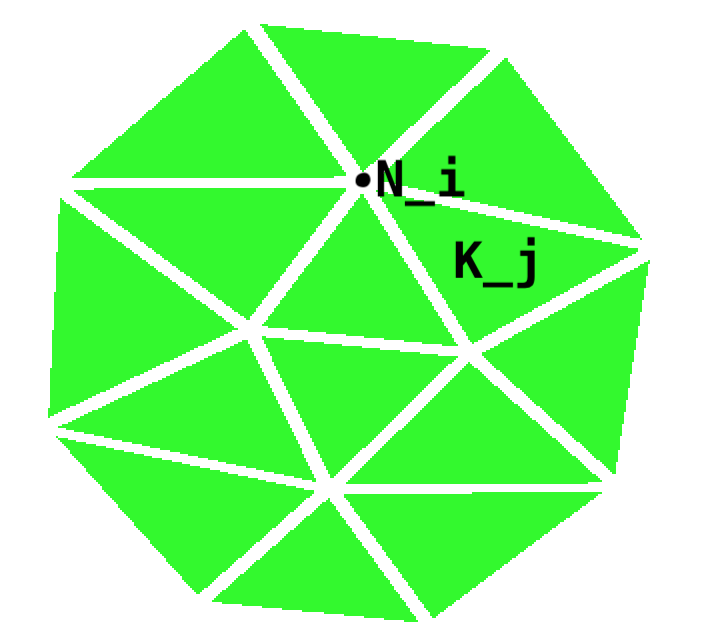}
    \caption{Triangulation of $\Omega$}
    \textit{Note:} Triangles have been separated in the edges to take a better look, but the triangulation has no empty spaces.
\end{figure}
The $h$ in the notation $T_h$ is important for the project because it gives a sense of the size for the mesh. It is defined as follows: $h=\max\{diam(K):K\in T_h\}$ where $diam(K)=\text{longest side of }K$.\\
\newline
Now, let $V_h=\{v:v \text{ continuous on }\Omega, v|_K \text{ linear for }K\in T_h,\ v=0\text{ on }\Gamma\}$.\\
If we consider the nodes ($N_1,\dots,N_M$) of the triangulation that are not on the boundary, since $v=0$ there, and define the functions $\varphi_j(N_i)=\left\{ \begin{array}{lcc}
             1 &   ,\ i=j\\
             \\0 &  ,\ i\not=j \\
             \end{array}
   \right.$
for $i,j=1,\dots,M$ in a way that $\varphi_j\in V_h$:
\begin{figure}[h!]
    \centering
    \includegraphics[scale=0.3]{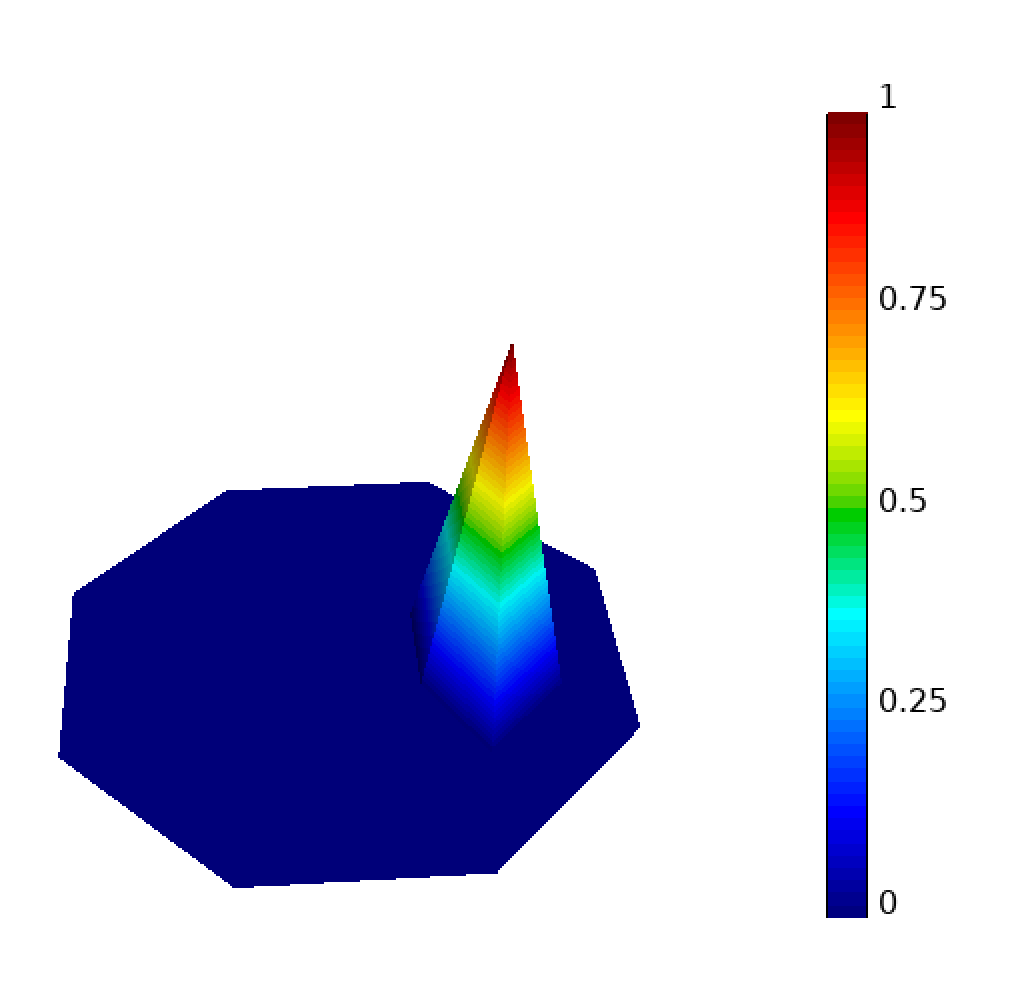}
    \includegraphics[scale=0.3]{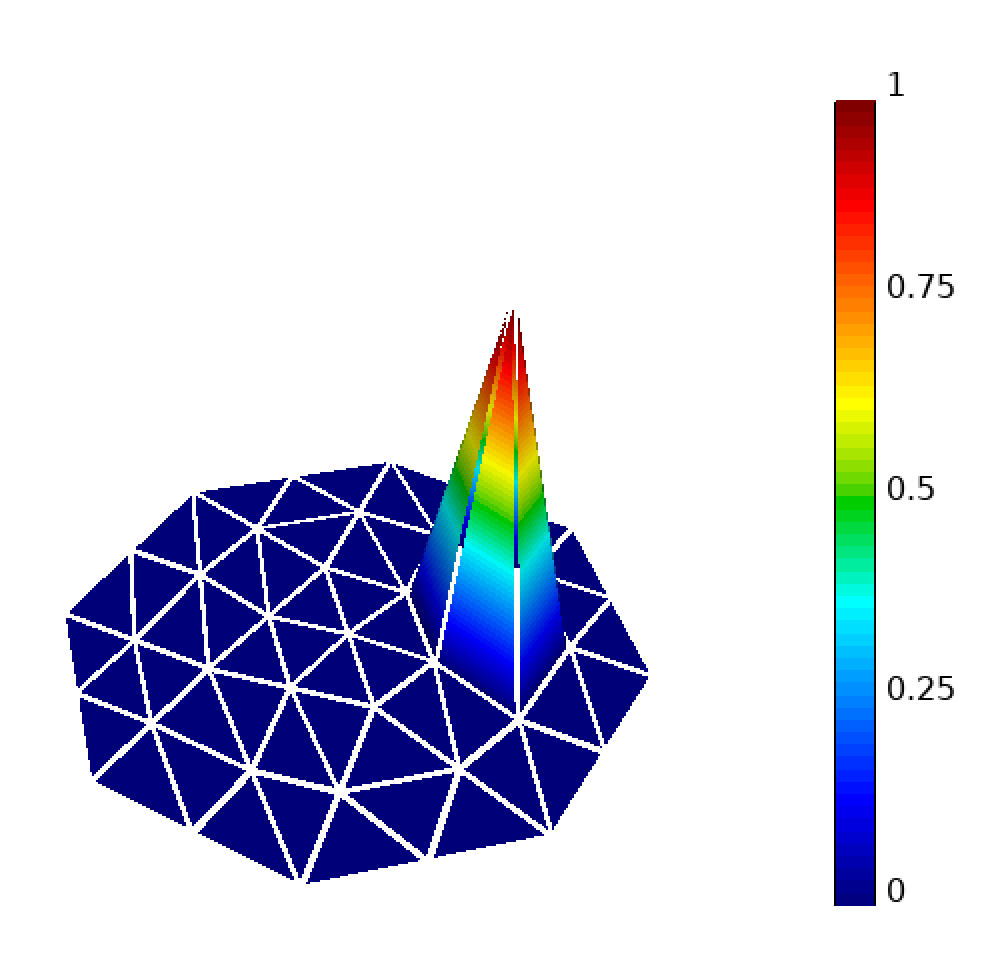}
    \caption{Function $\varphi_j$}
\end{figure}\\
With this, $V_h=gen\{\varphi_i:i=1,\dots,M\}$ because, for $v(x)\in V_h$, \\$v(x)=\sum_{j=1}^M\xi_j\varphi_j(x),$ with $\xi_j=v(N_j)\ and\ x\in\Omega\cup\Gamma$. So, $V_h$ is a finite-dimensional subspace of $V$. \cite{Lagrange}\\
\newline
Then, if $p_h\in V_h$ satisfies (4) for all $v\in V_h$ then, in particular:
\begin{equation}
    \int_\Omega\nabla p_h\cdot\nabla \varphi_j=\int_\Omega f\ \varphi_j,\ \ j=1,\dots,M
\end{equation}
As, $\nabla p_h=\sum_{i=1}^M\xi_i\nabla\varphi_i$ with $\xi_i=p_h(N_i)$, replacing on (5) we get:
\begin{equation}
    \sum_{i=1}^M\xi_i\int_\Omega\nabla\varphi_i\cdot\nabla\varphi_j=\int_\Omega f\ \varphi_j,\ \ j=1,\dots,M
\end{equation}
Finally, (6) is a linear system of $M$ equations and $M$ unknowns ($\xi_1,\dots,\xi_M$), which can be written as:
\begin{equation}
    A\xi=b
\end{equation}
where $A[i,j]=\int_\Omega\nabla\varphi_i\cdot\nabla\varphi_j$, $\xi[i]=p_h(N_i)$ and $b[i]=\int_\Omega f\ \varphi_i$.\\
\newline
In \cite{Lagrange}, it is shown that (7) has an unique solution and that matrix $A$ has useful properties for computing with it. Also, we can solve (7) with MFEM library.

\subsection{Mixed finite elements}
First, let's define some important spaces, where $\Omega$ is a bounded domain in $\mathbb{R}^2$ and $\Gamma$ its boundary \cite{Mixed}:
\begin{equation*}
    L^2(\Omega)=\{v:\Omega\rightarrow\mathbb{R}\ \Big|\int_\Omega v^2<\infty\}
\end{equation*}
\begin{equation*}
    H^1(\Omega)=\{v\in L^2(\Omega)\ \Big|\ \frac{\partial v}{\partial x},\frac{\partial v}{\partial y}\in L^2(\Omega)\}
\end{equation*}
\begin{equation*}
    H_0^1(\Omega)=\{v\in H^1(\Omega)\ |\ v=0\ on\ \Gamma\}
\end{equation*}
\begin{equation*}
    H(div;\Omega)=\{\mathbf{v}\in L^2(\Omega)\times L^2(\Omega)\ |\ div(\mathbf{v})\in L^2(\Omega)\}
\end{equation*}
As above, let $\Omega\in\mathbb{R}^2$ be a bounded domain with boundary $\Gamma$ and consider the following problem \cite{Mixed}:
\begin{equation}
\begin{split}
        -\Delta &p=f\text{ in } \Omega\\
        &p=0\text{ in }\Gamma
\end{split}
\end{equation}
where $f\in L^2(\Omega)$ and $\Delta p=\frac{\partial^2p}{\partial x^2}+\frac{\partial^2p}{\partial y^2}$.\\
\newline
This problem is the same problem considered in 2.1, but with a special condition for $f$, and can be reduced to:
\begin{equation*}
\int_\Omega\nabla v\cdot\nabla p=\int_\Omega f\ v,\text{ for all $v\in V$}
\end{equation*}
where Dirichlet boundary condition ($p=0\ in\ \Gamma$) is essential. \\
\newline
\texttt{Remark:} \textit{The space $V$ can be replaced with $H_0^1(\Omega)$ as seen in \cite{Mixed}}.\\
\newline
However, for mixed formulation, boundary won't be essential but natural:\\
\newline
Let $u=\nabla p$ in $\Omega$. \\
With this, problem (1) can be written as:
\begin{equation}
\begin{split}
    &u=\nabla p\text{ in }\Omega\\
    &div(u)=-f\text{ in } \Omega\\
        &p=0\text{ in }\Gamma
\end{split}
\end{equation}
because $\Delta p=div(\nabla p)$. \\
\newline
Now, following a similar procedure as in section 2.1:\\
\newline
Multiply the first equation of (2) by some $\mathbf{v}\in H(div;\Omega)$ and integrate both sides:
\begin{equation}
    \int_\Omega u\ \mathbf{v}=\int_\Omega\nabla p\cdot\mathbf{v}
\end{equation}
Consider Green's identity \cite{Mixed}:
\begin{equation}
    \int_\Omega \mathbf{v}\cdot\nabla p+\int_\Omega p\ div(\mathbf{v})=\int_\Gamma(\mathbf{v}\cdot\eta) p
\end{equation}
Replacing (4) in (3), and considering the third equation of (2), we get:
\begin{equation}
    \int_\Omega u\ \mathbf{v}+\int_\Omega p\ div(\mathbf{v})=\int_\Gamma(\mathbf{v}\cdot\eta) p
\end{equation}
where $\eta$ is the normal vector exterior to $\Gamma$.\\
\newline
On the other hand, we can multiply the second equation of (2) by some $w \in L^2(\Omega)$, integrate and obtain:
\begin{equation}
    \int_\Omega w\ div(u)=-\int_\Omega f\ w
\end{equation}
\textit{Remark: The boundary integral depends directly on the value of $p$ in $\Gamma$. And, this is what's called a \textbf{natural} boundary condition}.\\
\newline
Finally, applying boundary condition $p=0\ \text{ in }\Gamma$ into (5), and joining (5) and (6). We get the following problem deduced from (2):
\begin{equation}
    \begin{split}
        &\int_\Omega u\ \mathbf{v}+\int_\Omega p\ div(\mathbf{v})=0\\
        &\int_\Omega w\ div(u)=-\int_\Omega f\ w
    \end{split}
\end{equation}
\textbf{\textit{Note:}} For this problem, the objective is to find $(u,p)\in H(div;\Omega)\times L^2(\Omega)$ such that (7) is satisfied for all $\mathbf{v}\in H(div;\Omega),w\in L^2(\Omega)$.\\
\newline
For the discretized problem related to (7), define \cite{Mixed} the following spaces for a fixed \textit{triangulation} $T_h$ of the domain $\Omega$ and a fixed integer $k\geq0$:
\begin{equation*}
    \begin{split}
        &H_h:=\{\mathbf{v_h}\in H(div;\Omega):\mathbf{v_h}|_K\in RT_k(K)\text{ for all }K\in T_h\}\\
        &L_h:=\{w_h\in L^2(\Omega):w_h|_K\in \mathbb{P}_k(K)\text{ for all }K\in T_h\}
    \end{split}
\end{equation*}
where 
\begin{equation*}
    \begin{split}
        &\mathbb{P}_k(K)=\{p:K\rightarrow\mathbb{R}\ :\ p\text{ is a polynomial of degree }\leq k\}\\
        &RT_k(K)=[\mathbb{P}_k(K)\times\mathbb{P}_k(K)]+\mathbb{P}_k(K)x
    \end{split}
\end{equation*}
Note that $\mathbf{p}\in RT_k(K)$ if and only if there exist $p_0,p_1,p_2\in\mathbb{P}_k(K)$ such that
\begin{equation*}
    \mathbf{p}(x)=\begin{pmatrix}
    p_1(x)\\p_2(x)
    \end{pmatrix}
    +p_0(x)\begin{pmatrix}
    x\\y
    \end{pmatrix}
    \text{ for all }\begin{pmatrix}
    x\\y
    \end{pmatrix}\in K
\end{equation*}
Also, $\mathbf{p}$ has a degree of $k+1$.\\
\newline
Then, problem (7) can be changed to: find $(u_h,p_h)\in H_h\times L_h$ such that
\begin{equation}
    \begin{split}
        &\int_\Omega u_h\ \mathbf{v}_h+\int_\Omega p_h\ div(\mathbf{v}_h)=0\\
        &\int_\Omega w_h\ div(u_h)=-\int_\Omega f\ w_h
    \end{split}
\end{equation}
for all $\mathbf{v}_h\in H_h,w_h\in L_h$.\\
\newline
As spaces $H_h$ and $L_h$ are finite dimensional, they have a finite basis. That is, $H_h=gen\{\varphi_i:i=1,\dots,M\}$ and $L_h=gen\{\psi_j:j=1,\dots,N\}$. Then, $u_h=\sum_{i=i}^Mu_i\varphi_i$ and $p_h=\sum_{j=1}^Np_j\psi_j$, where $u_i$ and $p_j$ are scalars.\\
\newline
In particular, as $\varphi_k\in H_h$ and $\psi_l\in L_h$, we have that problem (8) can be written as
\begin{equation}
    \begin{split}
    &\int_\Omega\left(\sum_{i=i}^Mu_i\varphi_i\right)\varphi_k+\int_\Omega\left(\sum_{j=1}^Np_j\psi_j\right)div(\varphi_k)=0\\
    &\int_\Omega\psi_ldiv\left(\sum_{i=1}^Mu_i\varphi_i\right)=\int_\Omega f\psi_l\
    \end{split}
\end{equation}
for $k=1,\dots,M$ and $l=1,\dots,N$. Which is equivalent to the following by rearranging scalars:
\begin{equation}
    \begin{split}
    &\sum_{i=i}^Mu_i\int_\Omega\varphi_i\cdot\varphi_k+\sum_{j=1}^Np_j\int_\Omega\psi_jdiv(\varphi_k)=0\\
    &\sum_{i=i}^Mu_i\int_\Omega\psi_ldiv(\varphi_i)=\int_\Omega f\psi_l
    \end{split}
\end{equation}
for $k=1,\dots,M$ and $l=1,\dots,N$. This problem (10) can be formulated into the following matrix system
\begin{equation}
    \begin{pmatrix}
    A & B\\
    B^t & 0
    \end{pmatrix}\begin{pmatrix}
    U\\
    P
    \end{pmatrix}=\begin{pmatrix}
    0\\
    F
    \end{pmatrix}
\end{equation}
where $A$ is a $N\times N$ matrix, $B$ is a $M\times N$ matrix with $B^t$ it's transpose, $U$ is a $M$-dimensional column vector and $P,F$ are $N$-dimensional column vectors. \\
The entries of these arrays are $A[i,j]=\int_\Omega\varphi_i\cdot\varphi_j$, $B[i,j]=\int_\Omega\psi_jdiv(\varphi_i)$, $U[i]=u_i$, $P[i]=p_i$ and $F[i]=\int_\Omega f\psi_i$.\\
\newline 
(11) is a multilinear system that can be solved for $(U,P)$ with a computer using MFEM library. Note that with the entries of $U$ and $P$, the solution $(u_h,p_h)$ of (8) can be computed by their basis representation.\\
\newline
\textit{Note: The spaces defined to discretize the problem are called Raviart-Thomas finite element spaces. The fixed integer k is also called the order of the shape functions. And, the parameter $h$ is the same as in section 2.1, which is a meassure of size for $T_h$.} 
\subsection{Finite elements summary}
In sections 2.1 and 2.2 we studied two finite element methods. In general aspects, this is what was done:
\begin{itemize}
    \item Consider the problem of solving Poisson's equation with homogeneous Dirichlet boundary conditions. That is, the problem considered in previous sections.
    \item Multiply by some function (\textit{test function}) and integrate.
    \item Develop some equations applying boundary conditions.
    \item Discretize the domain.
    \item Define some finite-dimensional function spaces.
    \item Assemble the basis into the equation and form a matrix system.
\end{itemize}
The functions that form part of the finite-dimensional spaces are called $shape\ functions$. In Lagrange formulation, those where the functions in $V_h$, and in mixed formulation, those where the functions in $H_h$ and $L_h$.\\
\newline
The parameter $h$, denotes the size of the elements in the triangulation of the domain.\\
\newline
Both problems were solved with Dirichlet boundary condition ($= 0$). In Lagrange formulation it was essential, and in mixed formulation, it was natural.\\
\newline
In a more general aspect, the discretization of the space can be done without using triangles, but rather using quads or other figures.
\subsection{Higher order shape functions}
This is a very brief section that has the purpose of explaining a little bit of finite elements order, because in section 3 we will use different orders for the shape functions.\\
\newline
In general aspects, the order of a shape function is similar to the order of a polynomial. In mixed formulation we approached this when talking about Raviart-Thomas spaces, as in this spaces if the order of the polynomial is $k$, then the order of the shape function is $k+1$.\\
\newline
In the original introduction of the Lagrange formulation, the order of the shape functions was set to one. Better approximations can be obtained by using polynomials of higher order. Instead of defining 
\begin{equation*}
    V_h=\{v:v \text{ continuous on }\Omega, v|_K \text{ linear for }K\in T_h,\ v=0\text{ on }\Gamma\}
\end{equation*}
one can define, for a fixed order $k$: 
\begin{equation*}
    V^k_h=\{v:v \text{ continuous on }\Omega, v|_K \text{ polynomial of order at most }K\in T_h,\ v=0\text{ on }\Gamma\}.
\end{equation*}
\textit{Remark: For a fixed $k$, Lagrange shape functions have order 1 less than mixed shape functions.}\\
\newline
For example, as seen in \cite{Galvis}, the space of Bell triangular finite elements for a given triangulation $T_h$ is the space of functions that are polynomials of order 5 when restricted to every triangle $K\in T_h$. That is, if $v$ is in this space, then:
\begin{equation*}
    v|_K(x,y)=a_1x^5+a_2y^5+a_3x^4y+a_4xy^4+\dots+a_{16}x+a_{17}y+a_{18}
\end{equation*}
for all $K\in T_h$. Here, the constants $a_i,\ i=1,\dots,18$ correspond to $v$'s \textit{DOF} (degrees of freedom).
\begin{figure}[h!]
    \centering
    \includegraphics[scale=0.3]{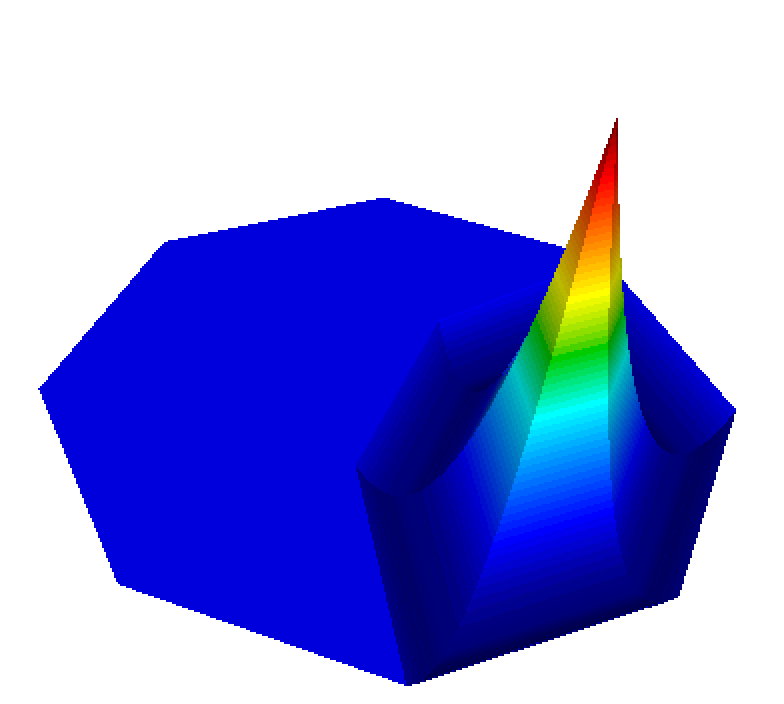}
    \includegraphics[scale=0.28]{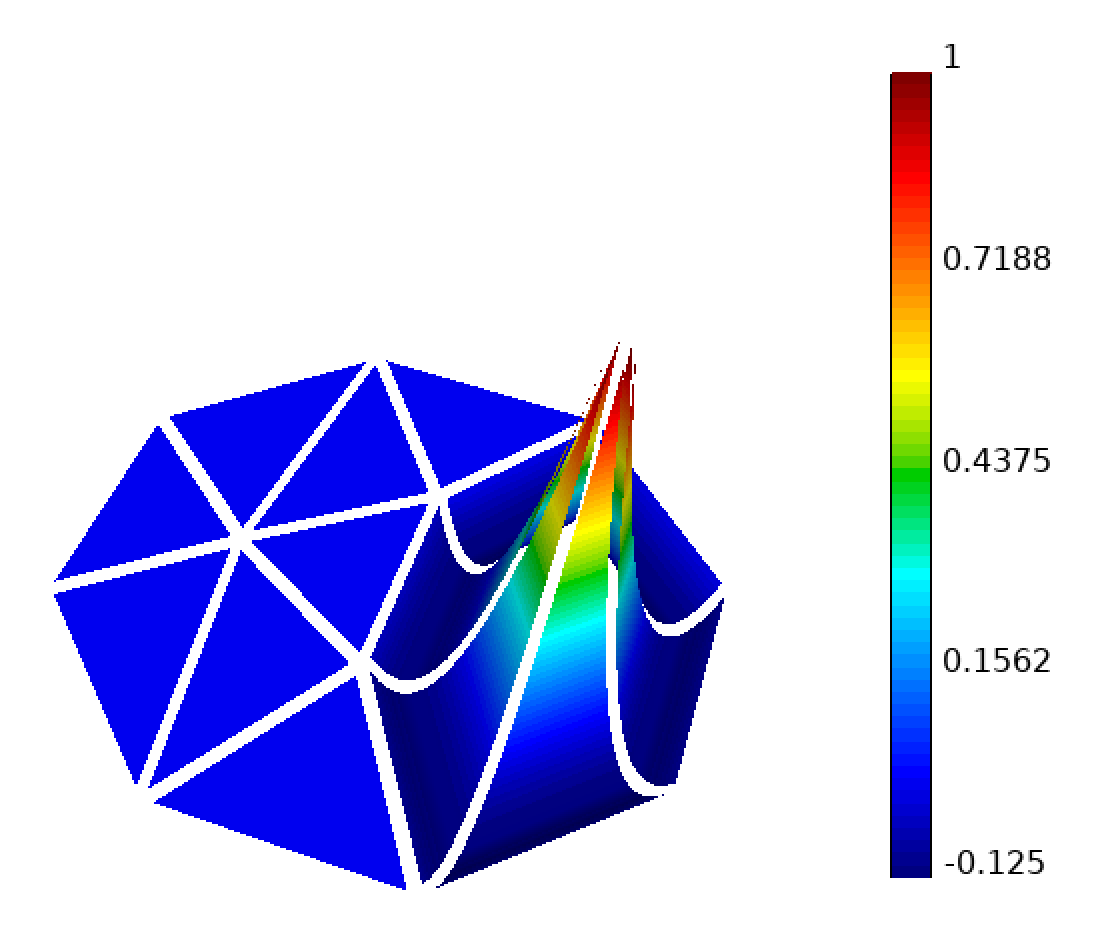}
    \caption{Finite element of order 2}
\end{figure}
\begin{figure}
    \centering
    \includegraphics[scale=0.35]{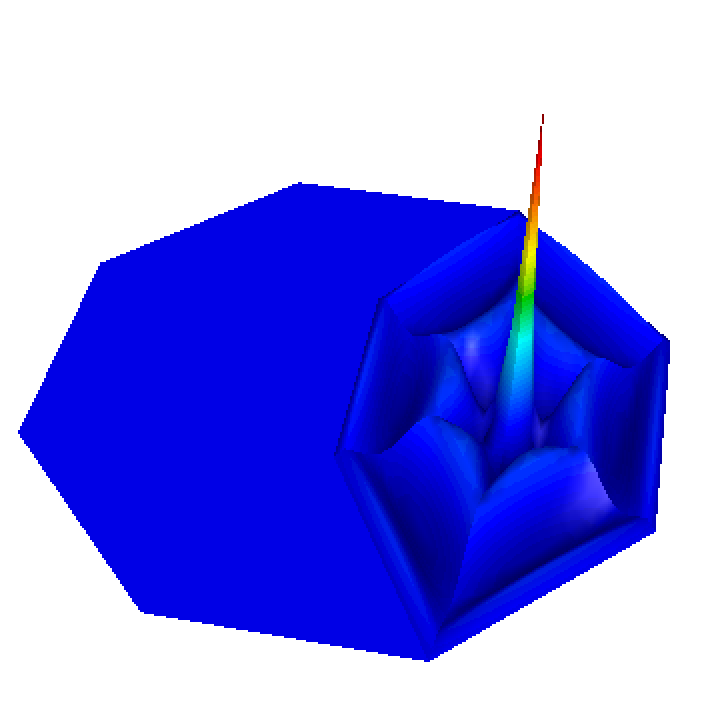}
    \includegraphics[scale=0.3]{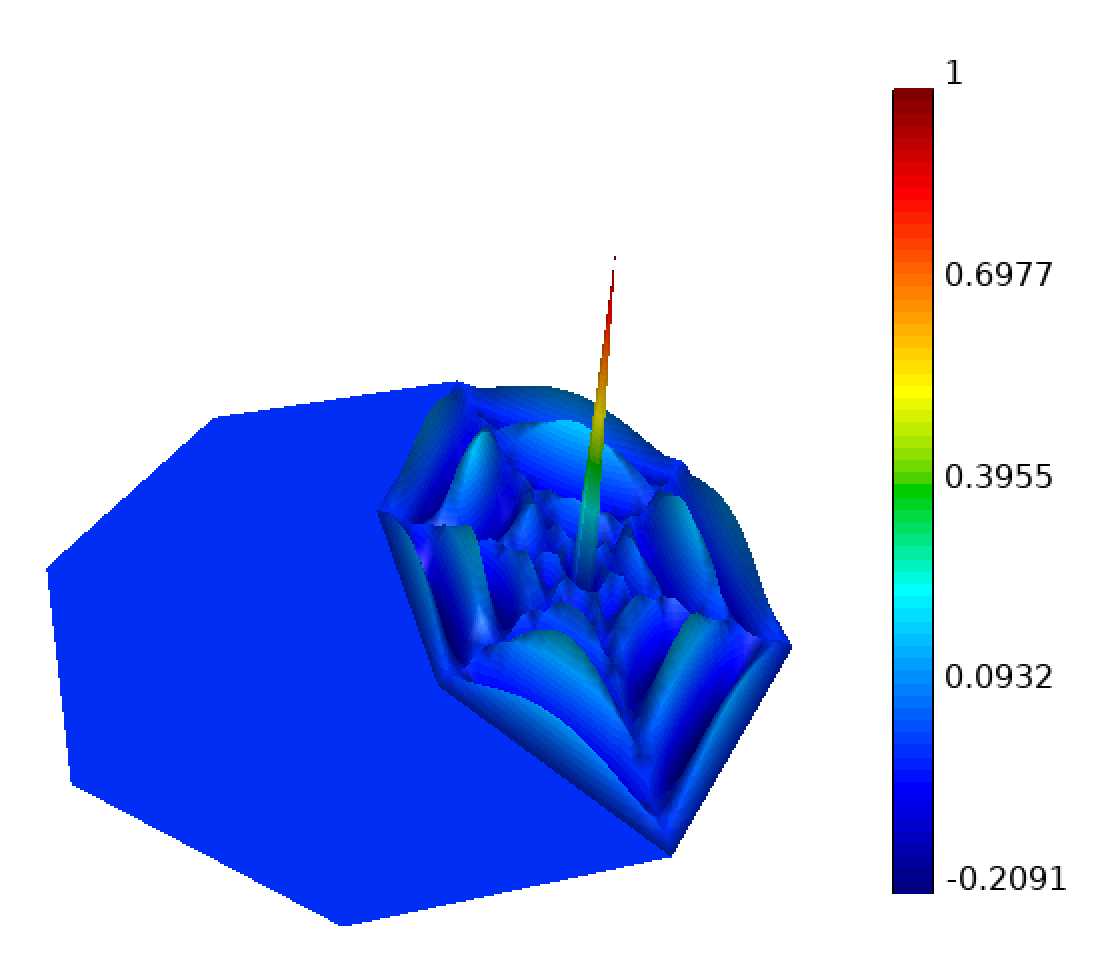}
    \caption{Finite elements of orders 5 (left) and 10 (right)}
\end{figure}
\newpage
\subsection{MFEM library}
In this project, we worked with MFEM's \textit{Example\#1} and \textit{Example\#5} which can be found on \cite{MFEM}. Example\#1 uses standard Lagrange finite elements and Example\#5 uses Raviart-Thomas mixed finite elements. Further, in section 3.1, we find the parameters so that both problems are equivalent and then (section 3.4), we compare the solutions.
\subsubsection{Overview}
According to it's official site \cite{MFEM}, MFEM is a free, lightweight, scalable C++ library for finite element methods that can work with arbitrary high-order finite element meshes and spaces.\\
\newline
MFEM has a serial version (which we are using) and a parallel version (for parallel computation). 
\newpage
The main classes (with a brief and superficial explanation of them) that we are going to use in the code are:
\begin{itemize}
    \item Mesh: domain with the partition.
    \item FiniteElementSpace: space of functions defined on the finite element mesh.
    \item GridFunction: mesh with values (solutions).
    \item $\_$Coefficient: values of GridFunctions or constants.
    \item LinearForm: maps an input function to a vector for the rhs.
    \item BilinearForm: used to create a global sparse finite element matrix for the lhs.
    \item $\_$Vector: vector.
    \item $\_$Solver: algorithm for solution calculation.
    \item $\_$Integrator: evaluates the bilinear form on element's level.
\end{itemize}
The ones that have $\_$ are various classes whose name ends up the same and work similarly.\\
\newline
\textit{\underline{Note:}} \\\textbf{lhs}: left hand side of the linear system.\\
\textbf{rhs}: right hand side of the linear system.

\subsubsection{Code structure}
An MFEM general code has the following steps (directly related classes with the step are written):
\begin{enumerate}
    \item Receive archive (.msh) input with the mesh and establish the order for the finite element spaces.
    \item Create mesh object, get the dimension, and refine the mesh (refinement is optional). \texttt{Mesh}
    \item Define the finite element spaces required. \texttt{FiniteElementSpace}
    \item Define coefficients, functions, and boundary conditions of the problem. \texttt{\textit{X}Coefficient}
    \item Define the LinearForm for the rhs and assemble it. \texttt{LinearForm, \textit{X}Integrator}
    \item Define the BilinearForm for the lhs and assemble it. \texttt{BilinearForm, \textit{X}Integrator}
    \item Solve the linear system. \texttt{\textit{X}Solver, \textit{X}Vector}
    \item Recover solution. \texttt{GridFunction}
    \item Show solution with a finite element visualization tool like \textbf{Glvis} (optional).
\end{enumerate}
\section{A case study}
In this section: we take examples 1 and 5 from \cite{MFEM}, define their problem parameters in such way that they're equivalent, create a code that implements both of them at the same time and compares both solutions ($L_2$ norm), run the code with different orders, and analyse the results.\\
\newline
Some considerations to have into account are:
\begin{itemize}
    \item For a fair comparison, order for Mixed method should be 1 less than order for Lagrange method. Because, with this, both shape functions would have the same degree.
    \item The code has more steps than shown in section 2.3.2 because we are running two methods and comparing solutions.
    \item We will compare pressures and velocities with respect to the order of the shape functions and the size of the mesh ($h$ parameter).
    \item For the problem, the exact solution is known, so, we will use it for comparison.
    \item The max order and refinement level to be tested is determined by our computational capacity (as long as solvers converge fast).
    \item The mesh used is a star with rhomboidal elements. 
\end{itemize}
\subsection{Problem}
Example\#1 \cite{MFEM}:
\begin{equation}
\begin{split}
        -\Delta &p=1\text{ in } \Omega\\
        &p=0\text{ in }\Gamma
\end{split}
\end{equation}
Example\#5 \cite{MFEM}:
\begin{equation}
\begin{split}
        &k\mathbf{u}+\nabla p=f\text{ in } \Omega\\
        &-div(\mathbf{u})=g\text{ in } \Omega\\
        &-p=p_0\text{ in }\Gamma
\end{split}
\end{equation}
From the first equation of (2):
\begin{equation}
    \mathbf{u}=\frac{f-\nabla p}{k}
\end{equation}
Then, replacing (3) on the second equation of (2):
\begin{equation}
    -div\left(\frac{f-\nabla p}{k}\right)=g
\end{equation}
If we set $k=1;\ f=0\ and\ g=-1$ in (4), we get:
\begin{equation}
    -\Delta p=1
\end{equation}
which is the first equation of (1). \\
\newline
So, setting \textbf{($*$)} $p_0=0,\ k=1;\ f=0\ and\ g=-1$ in (2), we get:
\begin{equation}
\begin{split}
        &\mathbf{u}+\nabla p=0\text{ in } \Omega\\
        &-div(\mathbf{u})=-1\text{ in } \Omega\\
        &-p=0\text{ in }\Gamma
\end{split}
\end{equation}
Notice that from the first equation we get that $\mathbf{u}=-\nabla p$. This is important because in problem (1) we don't get $\mathbf{u}$ solution from the method, so, in the code, we will have to find it from $p$'s derivatives.\\
\newline
In the code, we will set the value of the parameters in the way shown here, so that both problems are the same. As seen in (3)-(5), problem (6) is equivalent to problem (1) with the values assigned for coefficients and functions in \textbf{($*$)}.
\subsection{Code}
The first part of the code follows the structure mentioned in 2.3.2, but implemented for two methods at the same time (and with some extra lines for comparison purposes). Also, when defining boundary conditions, the \textit{essential} one is established different from the \textit{natural} one. And, after getting all the solutions, there's a second part of the code where solutions are compared between them and with the exact one. \\
\newline
\textit{Note:} \\The complete code with explanations can be found on the \texttt{Appendix A}.
\newpage
However, before taking a look into it, here's the convention used for important variable names along the code:\\
\newline
\textbf{\textit{Notation:}}
\begin{table}[h!]
\begin{tabular}{|c|c|}
\hline
\textbf{Variable Name} & \textbf{Object}                               \\ \hline
X\_space               & Finite element space X                        \\ \hline
X\_mixed               & Variable assigned to a mixed method related object    \\ \hline
u                      & Velocity solution                                     \\ \hline
p                      & Pressure solution                                    \\ \hline
X\_ex                  & Variable assigned to an exact solution object \\ \hline
\end{tabular}
\end{table}

\subsection{Tests}
The tests will be run on the following domain:
\begin{figure}[h!]
    \centering
    \includegraphics[scale=0.5]{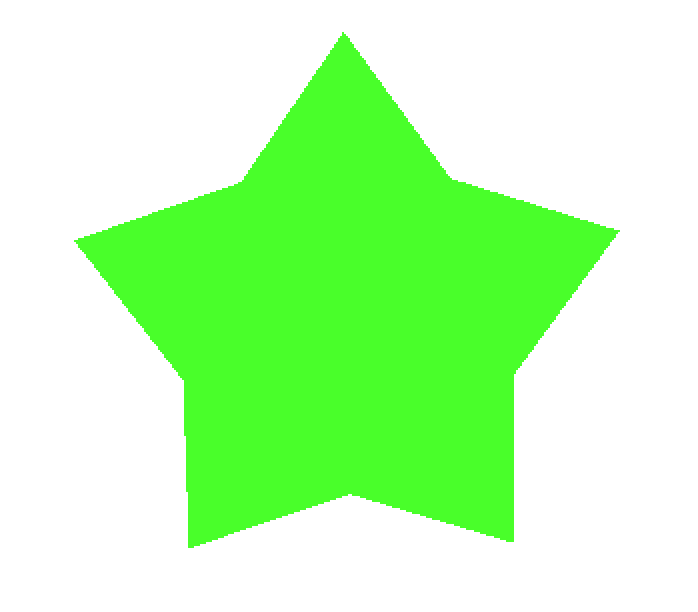}
    \caption{Star domain for tests}
\end{figure}\\
Each run test is determined by the \textit{order} of Lagrange shape functions and the \textit{h} parameter of the mesh. Remember that mixed shape functions have order equal to $\textit{order}-1$. The parameter \textit{order} is changed directly from the command line, while the parameter \textit{h} is changed via the number of times that the mesh is refined ($h=h(\#refinements)$). As we refine the mesh more times, finite elements of the partition decrease their size, and so, the parameter $h$ decreases.\\
\newline
Tests will be made with: $order=1,\dots,N$ and $refinements=0,\dots,M$, where $N,M$ depend on the computation capacity. The star mesh comes with a default partition which is shown below:
\newpage
\begin{figure}[h!]
    \centering
    \includegraphics[scale=0.5]{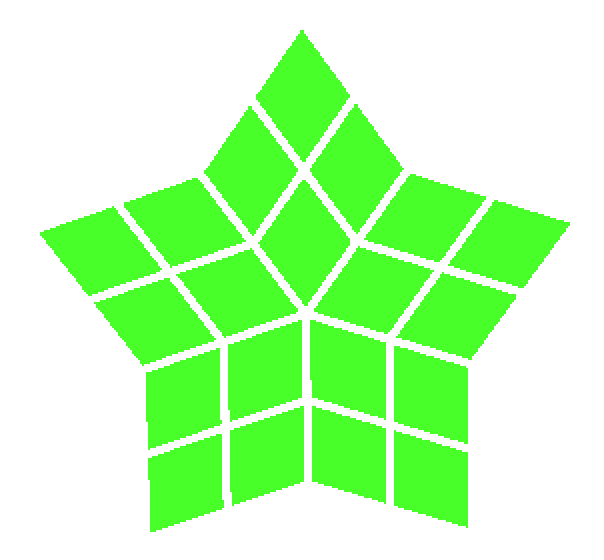}
    \caption{Mesh with no refinement}
\end{figure}
Results will be presented in graphs. However, all the exact values that were computed can be found in the \texttt{Appendix B}.
\subsection{Results}
Before showing the graphs, this is the output received in the visualization tool (Glvis) when running the code with $\textit{order}=2$ and $\#Refinements=3$ (graphically, Lagrange and Mixed solutions look the same):
\begin{figure}[h!]
    \centering
    \includegraphics[scale=0.25]{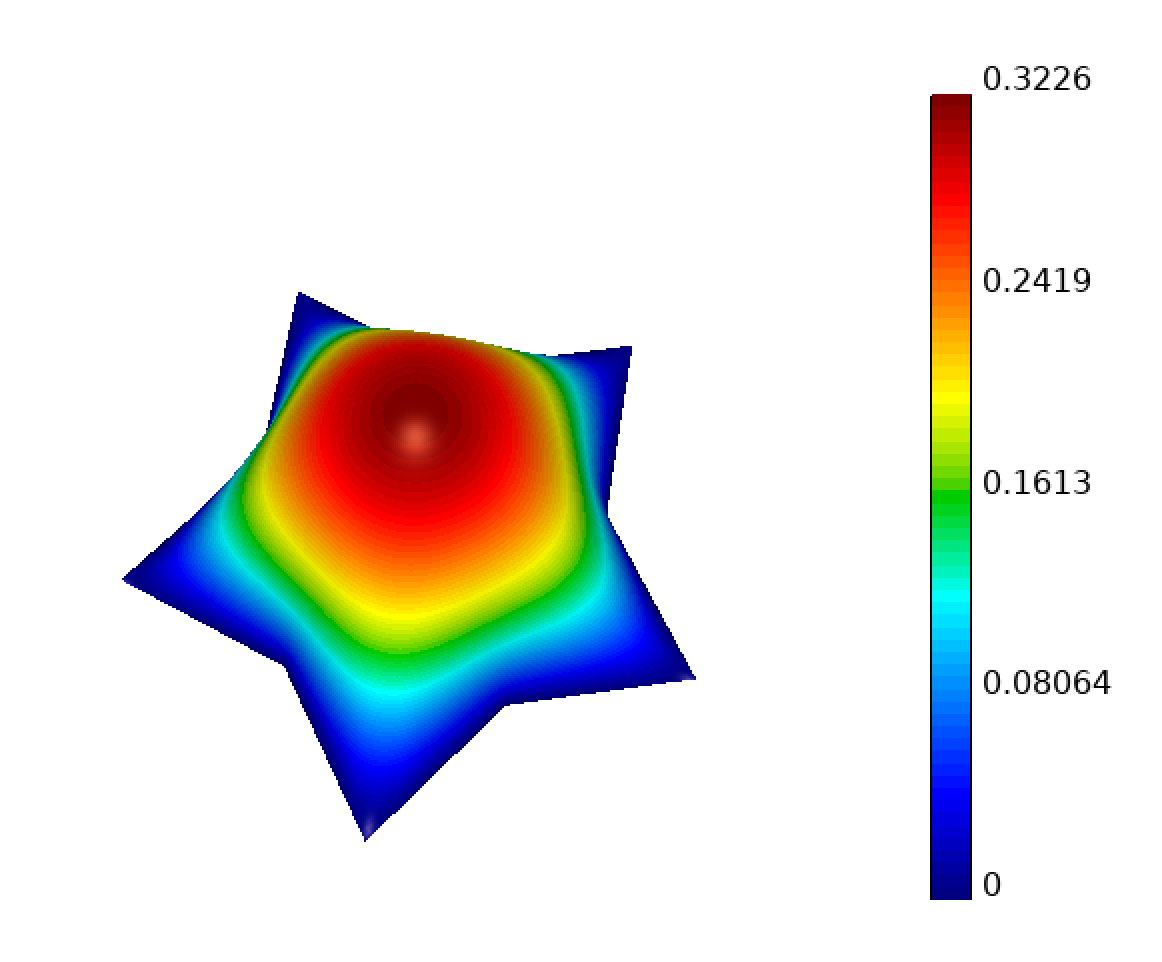}
    \includegraphics[scale=0.25]{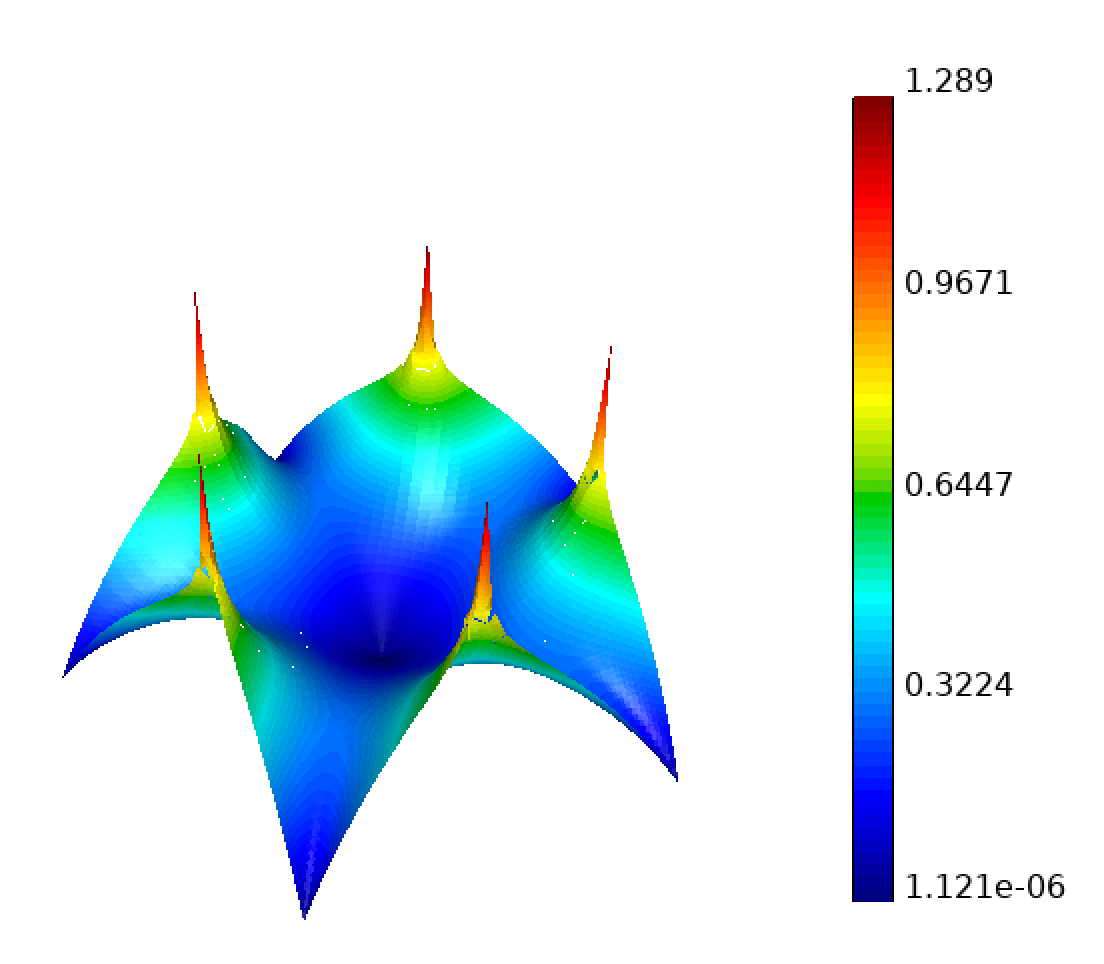}
    \caption{Glvis Visualization: Pressure (left) and Velocity (right)}
\end{figure}\\
\textit{Note:} Although velocity is a vector on each point, Glvis visualization tool doesn't shows it like that. It rather shows the $L^2$ norm of the vector.\\
\newline
In the following graphs, if $u=(u_x,u_y)$ is the solution obtained by the mixed or Lagrange finite element method and $u_{ex}=(u_{x_{ex}},u_{y_{ex}})$ is the exact solution for the problem, then:
\begin{equation*}
    U_{error} = \frac{\sqrt{\left(||u_x-u_{x_{ex}}||_{L^2}\right)^2+\left(||u_y-u_{y_{ex}}||_{L^2}\right)^2}}{||u_{ex}||_{L^2}}
\end{equation*}
\newpage
\begin{figure}[h!]
    \centering
    \includegraphics[scale = 0.5]{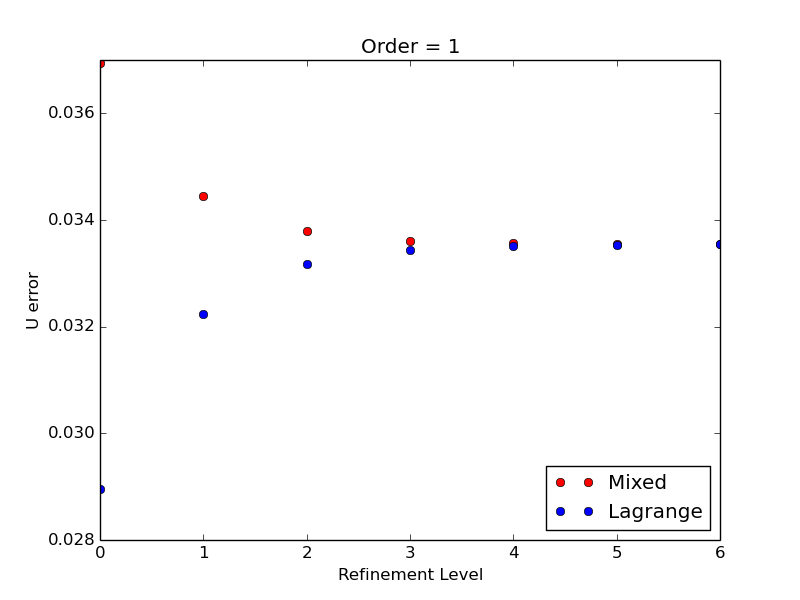}
    \caption{Order = 1}
\end{figure}
\begin{figure}[h!]
    \centering
    \includegraphics[scale = 0.5]{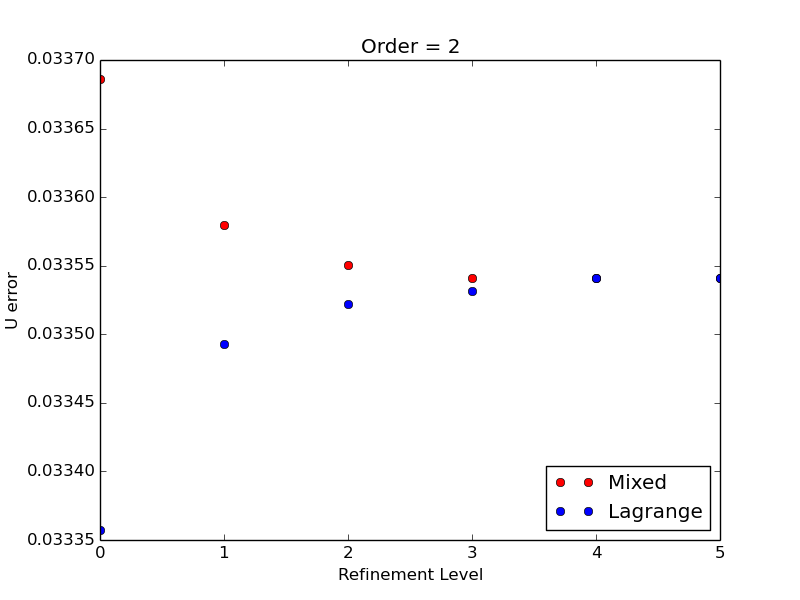}
    \caption{Order = 2}
\end{figure}
\newpage
\begin{figure}[h!]
    \centering
    \includegraphics[scale = 0.5]{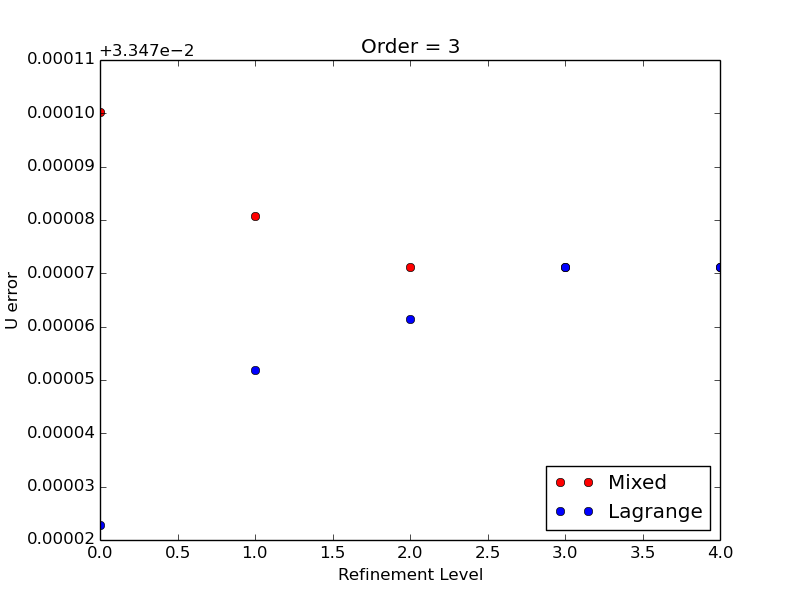}
    \caption{Order = 3}
\end{figure}
\begin{figure}[h!]
    \centering
    \includegraphics[scale = 0.5]{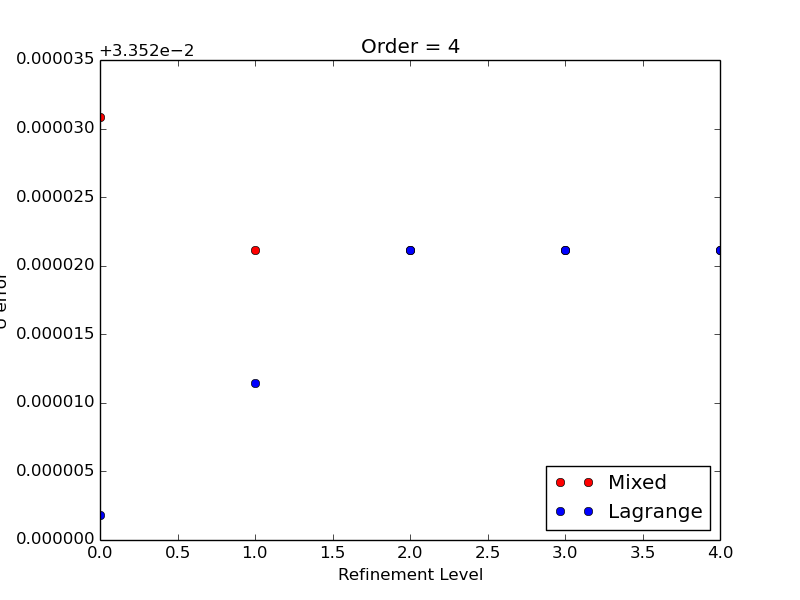}
    \caption{Order = 4}
\end{figure}
\newpage
\subsection{Analysis}
This section was done by analyzing the tables presented on the Appendix B.\\
\newline
To understand the information presented, take into account the following:
\begin{itemize}
    \item The exact solution would have value $1$ in \textit{X err}.
    \item If the two solutions obtained (Lagrange and Mixed) are exactly the same, the value in \textit{P comp} and \textit{U comp} would be $0$.
    \item Lower values of $h$ mean more mesh refinements, ie, smaller partition elements.
\end{itemize}
As it was expected, computational time increases as order and refinements increase.\\
\newline
Here are the most relevant observations that can be obtained after analysing the data corresponding to \textbf{absolute errors}:
\begin{itemize}
    \item For fixed order, absolute errors have little variation when reducing $h$ (max variation is $4.722$e$-03$ in $Uerr$ order 1).
    \item Absolute errors variation (respect to refinement) is lower when order is higher. For example; in order 2, $Perr$ is the same for each $h$ (up tu three decimal places); while in order 6, $Perr$ is the same for each $h$ (up to five decimal places).
    \item For fixed $h$, absolute errors remain almost constant between orders.
    \item $Perr$ (absolute error obtained for pressure with Lagrange) is always lower than $Pmx\ err$ (absolute error obtained for pressure with mixed).
    \item For fixed order, $Perr$ increases as $h$ decreases, while $Pmx\ err$ decreases as $h$ decreases.
    \item $Uerr$ (absolute error obtained for velocity with Lagrange) is always lower than $Umx\ err$ (absolute error obtained for velocity with mixed).
    \item For fixed order, $Uerr$ increases as $h$ decreases, while $Umx\ err$ decreases as $h$ decreases.
    \item As order increases, pressure absolute errors tend to be the same. In order 10, the difference between $Perr$ and $Pmx\ err$ is $0.000001$.
    \item As order increases, velocity absolute errors tend to be the same. In order 10, the difference between $Uerr$ and $Umx\ err$ is $<0.0000009$.
\end{itemize}

\newpage
And now, the most relevant observations that can be obtained after analysing the data corresponding to \textbf{comparison errors}:
\begin{itemize}
    \item Comparison errors, $Ucomp$ and $Pcomp$, decrease as $h$ decreases.
    \item When order increases, comparisons errors are lower for fixed $h$.
    \item Comparison error tends to $0$, as expected.
    \item Pressure comparison error lowers faster than velocity comparison error. Maximum comparison errors were found on order 1 with no refinements, where $Pcomp\approx7.5$e$-02$ and $Ucomp\approx3.7$e$-02$, and in minimum comparison errors were found on order 10 with 1 refinement (higher refinement level computed for order 10), where $Pcomp\approx5.1$e$-06$ and $Ucomp\approx9.8$e$-04$. It can be seen that $Pcomp$ improved in almost four decimal places while $Ucomp$ improved in just 2.
    \item For a fixed order, comparison error can be similar to a higher order comparison error, as long as enough refinements are made.
\end{itemize}

\section{Conclusion}
Adding up to the observations made in section 3.5, Lagrange solution and mixed solution tend to be the same when order and refinement levels increase, as expected. Also, Lagrange formulation is implemented more easily compared to mixed formulation but, with mixed formulation one can obtain pressure and velocity solutions at once. Furthermore, in MFEM, natural boundary conditions can be forced in an easier way compared to essential boundary conditions. Finally, it's important to note that finite element methods are a powerful mathematical tool used to solve potentially difficult problems.

\newpage
\section{Appendices}
\subsection{Appendix A}
Here, the code used (written in C++) is shown, with a brief explanations of it's functionality.\\
\lstset{language=C++,breaklines=true}
\xitem Include the required libraries (including MFEM) and begin main function.
\begin{lstlisting}[frame=single]
#include "mfem.hpp"
#include <fstream>
#include <iostream>
using namespace std;
using namespace mfem;
int main(int argc, char *argv[]){
\end{lstlisting}
\xitem Parse command-line options (in this project we only change "order" option) and print them.
\begin{lstlisting}[frame=single]
const char *mesh_file = "../data/star.mesh";
int order = 1;
bool visualization = true;
OptionsParser args(argc, argv);
args.AddOption(&mesh_file, "-m", "--mesh",
               "Mesh file to use.");
args.AddOption(&order, "-o", "--order",
            "Finite element order (polynomial degree).");
args.AddOption(&visualization, "-vis", "--visualization", "-no-vis", "--no-visualization",
               "Enable or disable GLVis visualization.");
args.Parse();
if (!args.Good()){
   args.PrintUsage(cout);
   return 1;
}
args.PrintOptions(cout);
\end{lstlisting}
\xitem Create mesh object from the star.mesh archive and get it's dimension.
\begin{lstlisting}[frame=single]
Mesh *mesh = new Mesh(mesh_file,1,1); 
int dim = mesh->Dimension();
\end{lstlisting}
\newpage
\xitem Refine the mesh a given number of times (uniform refinement).
\begin{lstlisting}[frame=single]
int ref_levels;
cout << "Refinements: ";
cin >> ref_levels;
for (int l = 0; l < ref_levels; l++){
mesh->UniformRefinement();
}
\end{lstlisting}
\xitem Get size indicator for mesh size (h\_max) and print it.
\begin{lstlisting}[frame=single]
double mesh_size, h = 0;
for (int i=0;i<mesh->GetNE();i++){
    mesh_size = mesh->GetElementSize(i,2);
    if(mesh_size>h){
        h = mesh_size;
    }
}
cout << "h: " << h << endl;
\end{lstlisting}
\xitem Define finite element spaces. For mixed finite element method, the order will be one less than for Lagrange finite element method. The last one is a vector L2 space that we will use later to get mixed velocity components.
\begin{lstlisting}[frame=single]
FiniteElementCollection *H1 = new H1_FECollection(order,dim);
FiniteElementSpace *H1_space = new FiniteElementSpace(mesh,H1);
FiniteElementCollection *hd(new RT_FECollection(order-1,dim));
FiniteElementCollection *l2(new L2_FECollection(order-1,dim));
FiniteElementSpace *Hdiv_space = new FiniteElementSpace(mesh,hd);
FiniteElementSpace *L2_space = new FiniteElementSpace(mesh,l2);
FiniteElementSpace *V_space = new FiniteElementSpace(mesh,l2,2);
\end{lstlisting}
\newpage
\xitem Define the parameters of the mixed problem. C functions are defined at the end. Boundary condition is natural.
\begin{lstlisting}[frame=single]
ConstantCoefficient k(1.0);
void fFun(const Vector & x, Vector & f);
VectorFunctionCoefficient fcoeff(dim, fFun);
double gFun(const Vector & x);
FunctionCoefficient gcoeff(gFun);
double f_bound(const Vector & x);
FunctionCoefficient fbndcoeff(f_bound);
\end{lstlisting}
\xitem Define the parameters of the Lagrange problem. Boundary condition is essential.
\begin{lstlisting}[frame=single]
ConstantCoefficient one(1.0);
Array<int> ess_tdof_list;
if (mesh->bdr_attributes.Size()){
   Array<int> ess_bdr(mesh->bdr_attributes.Max());
   ess_bdr = 1;
   H1_space->GetEssentialTrueDofs(ess_bdr, ess_tdof_list);
}
\end{lstlisting}
\xitem Define the exact solution. C functions are defined at the end.
\begin{lstlisting}[frame=single]
void u_ex(const Vector & x, Vector & u);
double p_ex(const Vector & x);
double u_ex_x(const Vector & x);
double u_ex_y(const Vector & x);
\end{lstlisting}
\xitem Get space dimensions and crate vectors for the right hand side.
\begin{lstlisting}[frame=single]
Array<int> block_offsets(3);
block_offsets[0] = 0;
block_offsets[1] = Hdiv_space->GetVSize();
block_offsets[2] = L2_space->GetVSize();
block_offsets.PartialSum();
BlockVector rhs_mixed(block_offsets);
Vector rhs(H1_space->GetVSize());
\end{lstlisting}
\newpage
\xitem Define the right hand side. These are LinearForm objects associated to some finite element space and rhs vector. "f" and "g" are for the mixed method and "b" is for the other method. "rhs" vectors are the variables that store the information of the right hand side.
\begin{lstlisting}[frame=single]
LinearForm *fform(new LinearForm);
fform->Update(Hdiv_space, rhs_mixed.GetBlock(0), 0);
fform->AddDomainIntegrator(new VectorFEDomainLFIntegrator(fcoeff));
fform->AddBoundaryIntegrator(new VectorFEBoundaryFluxLFIntegrator(fbndcoeff));
fform->Assemble();

LinearForm *gform(new LinearForm);
gform->Update(L2_space, rhs_mixed.GetBlock(1), 0);
gform->AddDomainIntegrator(new DomainLFIntegrator(gcoeff));
gform->Assemble();

LinearForm *b(new LinearForm);
b->Update(H1_space, rhs, 0);
b->AddDomainIntegrator(new DomainLFIntegrator(one));
b->Assemble();
\end{lstlisting}
\xitem Create variables to store the solution. "x" is the vector used as input in the iterative method.
\begin{lstlisting}[frame=single]
BlockVector x_mixed(block_offsets);
GridFunction u_mixed(Hdiv_space), p_mixed(L2_space), ux_mixed(L2_space), uy_mixed(L2_space), ue(V_space);
Vector x(H1_space->GetVSize());
GridFunction ux(L2_space),uy(L2_space),p(H1_space);
\end{lstlisting}
\newpage
\xitem Define the left hand side for mixed method. This is the bilinear form representing the Darcy matrix. VectorFEMMassIntegrator is asociated to $k*u-\nabla p$ and VectorFEDDivergenceIntegrator is asociated to $div(u)$.
\begin{lstlisting}[frame=single]
BilinearForm *mVarf(new BilinearForm(Hdiv_space));
MixedBilinearForm *bVarf(new MixedBilinearForm(Hdiv_space, L2_space));
mVarf->AddDomainIntegrator(new VectorFEMassIntegrator(k));
mVarf->Assemble();
mVarf->Finalize();
SparseMatrix &M(mVarf->SpMat());
bVarf->AddDomainIntegrator(new VectorFEDivergenceIntegrator);
bVarf->Assemble();
bVarf->Finalize();
SparseMatrix & B(bVarf->SpMat());
B *= -1.;
SparseMatrix *BT = Transpose(B);
BlockMatrix D(block_offsets);
D.SetBlock(0,0, &M);
D.SetBlock(0,1, BT);
D.SetBlock(1,0, &B);
\end{lstlisting}
\xitem Define the left hand side for Lagrange method. This is the bilinear form asociated to the laplacian operator. DiffusionIntegrator is asociated to $\Delta u$. The method FormLinearSystem is only used to establish the essential boundary condition.
\begin{lstlisting}[frame=single]
OperatorPtr A;
Vector XX,BB;
BilinearForm *a(new BilinearForm(H1_space));
a->AddDomainIntegrator(new DiffusionIntegrator(one));
a->Assemble();
a->FormLinearSystem(ess_tdof_list, p, *b, A, XX, BB);
\end{lstlisting}
\newpage
\xitem Solve linear systems with MINRES (for mixed) and CG (for Lagrange). SetOperator method establishes the lhs. Mult method executes the iterative algorithm and receives as input: the rhs and the vector to store the solution. Then convergence result is printed.
\begin{lstlisting}[frame=single]    
int maxIter(10000); 
double rtol(1.e-6);
double atol(1.e-10);
    
MINRESSolver Msolver;
Msolver.SetAbsTol(atol);
Msolver.SetRelTol(rtol);
Msolver.SetMaxIter(maxIter);
Msolver.SetPrintLevel(0); 
Msolver.SetOperator(D); 
x_mixed = 0.0;
Msolver.Mult(rhs_mixed, x_mixed);
if (Msolver.GetConverged())
   std::cout << "MINRES converged in " << Msolver.GetNumIterations() << " iterations with a residual norm of " << Msolver.GetFinalNorm() << ".\n";
else
   std::cout << "MINRES did not converge in " << Msolver.GetNumIterations() << " iterations. Residual norm is " << Msolver.GetFinalNorm() << ".\n";
   
CGSolver Lsolver;
Lsolver.SetAbsTol(atol);
Lsolver.SetRelTol(rtol);
Lsolver.SetMaxIter(maxIter);
Lsolver.SetPrintLevel(0);
Lsolver.SetOperator(*A);
x = 0.0;
Lsolver.Mult(rhs,x);
if (Lsolver.GetConverged())
   std::cout << "CG converged in " << Lsolver.GetNumIterations() << " iterations with a residual norm of " << Lsolver.GetFinalNorm() << ".\n";
else
   std::cout << "CG did not converge in " << Lsolver.GetNumIterations() << " iterations. Residual norm is " << Lsolver.GetFinalNorm() << ".\n";
\end{lstlisting}
\newpage
\xitem Save the solution into GridFunctions, which are used for error computation and visualization.
\begin{lstlisting}[frame=single]  
u_mixed.MakeRef(Hdiv_space, x_mixed.GetBlock(0), 0);
p_mixed.MakeRef(L2_space, x_mixed.GetBlock(1), 0);
p.MakeRef(H1_space,x,0);
\end{lstlisting}
\xitem Get missing velocities from the solutions obtained. Remember that $u = -\nabla p$. Mixed components are extracted using the auxiliary variable "ue" defined before.
\begin{lstlisting}[frame=single]  
p.GetDerivative(1,0,ux); 
p.GetDerivative(1,1,uy); 
ux *= -1;
uy *= -1;

VectorGridFunctionCoefficient uc(&u_mixed);
ue.ProjectCoefficient(uc);
GridFunctionCoefficient ux_mixed_coeff(&ue,1);
GridFunctionCoefficient uy_mixed_coeff(&ue,2);
ux_mixed.ProjectCoefficient(ux_mixed_coeff);
uy_mixed.ProjectCoefficient(uy_mixed_coeff);
\end{lstlisting}
\xitem Create the asociated Coefficient objects for error computation.
\begin{lstlisting}[frame=single]  
GridFunction* pp = &p;
GridFunctionCoefficient p_coeff(pp);
GridFunction* uxp = &ux;
GridFunction* uyp = &uy;
GridFunctionCoefficient ux_coeff(uxp);
GridFunctionCoefficient uy_coeff(uyp);
FunctionCoefficient pex_coeff(p_ex); 
VectorFunctionCoefficient uex_coeff(dim,u_ex); 
FunctionCoefficient uex_x_coeff(u_ex_x);
FunctionCoefficient uex_y_coeff(u_ex_y);
\end{lstlisting}
\xitem Define integration rule.
\begin{lstlisting}[frame=single]  
int order_quad = max(2, 2*order+1);
const IntegrationRule *irs[Geometry::NumGeom];
for (int i=0; i < Geometry::NumGeom; ++i){
    irs[i] = &(IntRules.Get(i, order_quad));
}
\end{lstlisting}
\newpage
\xitem Compute exact solution norms.
\begin{lstlisting}[frame=single]  
double norm_p = ComputeLpNorm(2., pex_coeff, *mesh, irs);
double norm_u = ComputeLpNorm(2., uex_coeff, *mesh, irs); 
double norm_ux = ComputeLpNorm(2., uex_x_coeff, *mesh, irs); 
double norm_uy = ComputeLpNorm(2., uex_y_coeff, *mesh, irs); 
\end{lstlisting}
\xitem Compute absolute errors and print them.
\begin{lstlisting}[frame=single]  
double abs_err_u_mixed = u_mixed.ComputeL2Error(uex_coeff,irs); 
printf("Velocity Mixed Absolute Error: %e\n", abs_err_u_mixed / norm_u);
double abs_err_p_mixed = p_mixed.ComputeL2Error(pex_coeff,irs); 
printf("Pressure Mixed Absolute Error: %e\n", abs_err_p_mixed / norm_p);
double abs_err_p = p.ComputeL2Error(pex_coeff,irs);
printf("Pressure Absolute Error: %e\n", abs_err_p / norm_p);
double abs_err_ux = ux.ComputeL2Error(uex_x_coeff,irs);
double abs_err_uy = uy.ComputeL2Error(uex_y_coeff,irs); 
double abs_err_u = pow(pow(abs_err_ux,2)+pow(abs_err_uy,2),0.5); 
printf("Velocity Absolute Error: %e\n", abs_err_u / norm_u);
\end{lstlisting}
\xitem Compute and print comparison errors.
\begin{lstlisting}[frame=single]  
double err_ux = ux_mixed.ComputeL2Error(ux_coeff,irs);
double err_uy = uy_mixed.ComputeL2Error(uy_coeff,irs);
double err_u = pow(pow(err_ux,2)+pow(err_uy,2),0.5);
printf("Velocity Comparison Error: %e\n", err_u / norm_u);
double err_p = p_mixed.ComputeL2Error(p_coeff, irs);
printf("Pressure Comparison Error: %e\n", err_p / norm_p);
\end{lstlisting}
\newpage
\xitem Visualize the solutions and the domain.
\begin{lstlisting}[frame=single]  
char vishost[] = "localhost";
int  visport   = 19916;
if(visualization){
    Vector x_domain(H1_space->GetVSize());
    GridFunction domain(H1_space);
    x_domain=0.0;
    domain.MakeRef(H1_space,x_domain,0);
    socketstream dom_sock(vishost, visport);
    dom_sock.precision(8);
    dom_sock << "solution\n" << *mesh << domain << "window_title 'Domain'" << endl;
    
    socketstream um_sock(vishost, visport);
    um_sock.precision(8);
    um_sock << "solution\n" << *mesh << u_mixed << "window_title 'Velocity Mixed'" << endl;
    socketstream pm_sock(vishost, visport);
    pm_sock.precision(8);
    pm_sock << "solution\n" << *mesh << p_mixed << "window_title 'Pressure Mixed'" << endl;
    socketstream uxm_sock(vishost, visport);
    uxm_sock.precision(8);
    uxm_sock << "solution\n" << *mesh << ux_mixed << "window_title 'X Velocity Mixed'" << endl;
    socketstream uym_sock(vishost, visport);
    uym_sock.precision(8);
    uym_sock << "solution\n" << *mesh << uy_mixed << "window_title 'Y Velocity Mixed'" << endl;
    
    socketstream p_sock(vishost, visport);
    p_sock.precision(8);
    p_sock << "solution\n" << *mesh << p << "window_title 'Pressure'" << endl;
    socketstream ux_sock(vishost, visport);
    ux_sock.precision(8);
    ux_sock << "solution\n" << *mesh << ux << "window_title 'X Velocity'" << endl;
    socketstream uy_sock(vishost, visport);
    uy_sock.precision(8);
    uy_sock << "solution\n" << *mesh << uy << "window_title 'Y Velocity'" << endl;
}
}
\end{lstlisting}
\xitem Define C functions.
\begin{lstlisting}[frame=single]  
void fFun(const Vector & x, Vector & f){
    f = 0.0;
}
double gFun(const Vector & x){
    return -1.0;
}
double f_bound(const Vector & x){
    return 0.0;
}
void u_ex(const Vector & x, Vector & u){
   double xi(x(0));
   double yi(x(1));
   double zi(0.0);
   u(0) = - exp(xi)*sin(yi)*cos(zi);
   u(1) = - exp(xi)*cos(yi)*cos(zi);
}
double u_ex_x(const Vector & x){
   double xi(x(0));
   double yi(x(1));
   double zi(0.0);
   return -exp(xi)*sin(yi)*cos(zi);
}
double u_ex_y(const Vector & x){
   double xi(x(0));
   double yi(x(1));
   double zi(0.0);
   return -exp(xi)*cos(yi)*cos(zi);
}
double p_ex(const Vector & x){
   double xi(x(0));
   double yi(x(1));
   double zi(0.0);
   return exp(xi)*sin(yi)*cos(zi);
}
\end{lstlisting}

\newpage
\subsection{Appendix B}
The \textit{order} parameter will be fixed for each table and $h$ parameter is shown in the first column. To interpret the results take into account that \textbf{P} refers to pressure, \textbf{U} refers to velocity, \textbf{mx} refers to mixed (from mixed finite element method), \textbf{err} refers to absolute error (compared to the exact solution), and \textbf{comp} refers to comparison (the error between the two solutions obtained by the two different methods).\\
\newline
\indent
\textbf{\textit{Order = 1}}
\begin{table}[h!]
\scalebox{0.9}{
\begin{tabular}{|c|c|c|c|c|c|c|}
\hline
\textbf{h} & \textbf{P comp} & \textbf{P err} & \textbf{Pmx err} & \textbf{U comp} & \textbf{U err} & \textbf{U mx err} \\ \hline
0.572063   & 7.549479e-02    & 1.021287e+00   & 1.025477e+00     & 3.680827e-02    & 1.029378e+00   & 1.037635e+00      \\ \hline
0.286032   & 3.627089e-02    & 1.022781e+00   & 1.023990e+00     & 1.727281e-02    & 1.032760e+00   & 1.035055e+00      \\ \hline
0.143016   & 1.791509e-02    & 1.023236e+00   & 1.023596e+00     & 9.222996e-03    & 1.033725e+00   & 1.034369e+00      \\ \hline
0.0715079  & 8.922939e-03    & 1.023372e+00   & 1.023480e+00     & 5.111295e-03    & 1.033999e+00   & 1.034182e+00      \\ \hline
0.035754   & 4.455715e-03    & 1.023412e+00   & 1.023445e+00     & 2.859769e-03    & 1.034077e+00   & 1.034130e+00      \\ \hline
0.017877   & 2.226845e-03    & 1.023424e+00   & 1.023435e+00     & 1.603788e-03    & 1.034100e+00   & 1.034115e+00      \\ \hline
\end{tabular}}
\end{table}

\textbf{\textit{Order = 2}}
\begin{table}[h!]
\scalebox{0.9}{
\begin{tabular}{|c|c|c|c|c|c|c|}
\hline
\textbf{h} & \textbf{P comp} & \textbf{P err} & \textbf{Pmx err} & \textbf{U comp} & \textbf{U err} & \textbf{U mx err} \\ \hline
0.572063   & 8.069013e-03    & 1.023329e+00   & 1.023554e+00     & 1.399079e-02    & 1.033924e+00   & 1.034255e+00      \\ \hline
0.286032   & 2.138257e-03    & 1.023391e+00   & 1.023470e+00     & 7.845012e-03    & 1.034056e+00   & 1.034146e+00      \\ \hline
0.143016   & 5.704347e-04    & 1.023417e+00   & 1.023442e+00     & 4.400448e-03    & 1.034093e+00   & 1.034120e+00      \\ \hline
0.0715079  & 1.537926e-04    & 1.023426e+00   & 1.023434e+00     & 2.469526e-03    & 1.034104e+00   & 1.034112e+00      \\ \hline
0.035754   & 4.194302e-05    & 1.023428e+00   & 1.023431e+00     & 1.385966e-03    & 1.034107e+00   & 1.034110e+00      \\ \hline
\end{tabular}}
\end{table}

\textbf{\textit{Order = 3}}
\begin{table}[h!]
\scalebox{0.9}{
\begin{tabular}{|c|c|c|c|c|c|c|}
\hline
\textbf{h} & \textbf{P comp} & \textbf{P err} & \textbf{Pmx err} & \textbf{U comp} & \textbf{U err} & \textbf{U mx err} \\ \hline
0.572063   & 8.691241e-04    & 1.023389e+00   & 1.023471e+00     & 8.745151e-03    & 1.034060e+00   & 1.034143e+00      \\ \hline
0.286032   & 2.477673e-04    & 1.023417e+00   & 1.023443e+00     & 4.911967e-03    & 1.034094e+00   & 1.034120e+00      \\ \hline
0.143016   & 7.316263e-05    & 1.023426e+00   & 1.023434e+00     & 2.756849e-03    & 1.034104e+00   & 1.034112e+00      \\ \hline
0.0715079  & 2.178864e-05    & 1.023428e+00   & 1.023431e+00     & 1.547232e-03    & 1.034108e+00   & 1.034110e+00      \\ \hline
\end{tabular}}
\end{table}

\textbf{\textit{Order = 4}}
\begin{table}[h!]
\scalebox{0.9}{
\begin{tabular}{|c|c|c|c|c|c|c|}
\hline
\textbf{h} & \textbf{P comp} & \textbf{P err} & \textbf{Pmx err} & \textbf{U comp} & \textbf{U err} & \textbf{U mx err} \\ \hline
0.572063   & 3.199774e-04    & 1.023412e+00   & 1.023448e+00     & 6.119857e-03    & 1.034088e+00   & 1.034124e+00      \\ \hline
0.286032   & 9.547574e-05    & 1.023424e+00   & 1.023435e+00     & 3.434952e-03    & 1.034103e+00   & 1.034114e+00      \\ \hline
0.143016   & 2.862666e-05    & 1.023428e+00   & 1.023431e+00     & 1.927814e-03    & 1.034107e+00   & 1.034111e+00      \\ \hline
\end{tabular}}
\end{table}

\newpage
\textbf{\textit{Order = 5}}
\begin{table}[h!]
\scalebox{0.9}{
\begin{tabular}{|c|c|c|c|c|c|c|}
\hline
\textbf{h} & \textbf{P comp} & \textbf{P err} & \textbf{Pmx err} & \textbf{U comp} & \textbf{U err} & \textbf{U mx err} \\ \hline
0.572063   & 1.552006e-04    & 1.023420e+00   & 1.023439e+00     & 4.578518e-03    & 1.034099e+00   & 1.034117e+00      \\ \hline
0.286032   & 4.658038e-05    & 1.023427e+00   & 1.023433e+00     & 2.569749e-03    & 1.034106e+00   & 1.034112e+00      \\ \hline
0.143016   & 1.406993e-05    & 1.023429e+00   & 1.023431e+00     & 1.442205e-03    & 1.034108e+00   & 1.034110e+00      \\ \hline
\end{tabular}}
\end{table}

\textbf{\textit{Order = 6}}
\begin{table}[h!]
\scalebox{0.9}{
\begin{tabular}{|c|c|c|c|c|c|c|}
\hline
\textbf{h} & \textbf{P comp} & \textbf{P err} & \textbf{Pmx err} & \textbf{U comp} & \textbf{U err} & \textbf{U mx err} \\ \hline
0.572063   & 8.612580e-05    & 1.023424e+00   & 1.023435e+00     & 3.584133e-03    & 1.034103e+00   & 1.034114e+00      \\ \hline
0.286032   & 2.600417e-05    & 1.023428e+00   & 1.023431e+00     & 2.011608e-03    & 1.034107e+00   & 1.034111e+00      \\ \hline
0.143016   & 7.897631e-06    & 1.023429e+00   & 1.023430e+00     & 1.128989e-03    & 1.034109e+00   & 1.034110e+00      \\ \hline
\end{tabular}}
\end{table}

\textbf{\textit{Order = 7}}
\begin{table}[h!]
\scalebox{0.9}{
\begin{tabular}{|c|c|c|c|c|c|c|}
\hline
\textbf{h} & \textbf{P comp} & \textbf{P err} & \textbf{Pmx err} & \textbf{U comp} & \textbf{U err} & \textbf{U mx err} \\ \hline
0.572063   & 5.243187e-05    & 1.023426e+00   & 1.023433e+00     & 2.899307e-03    & 1.034105e+00   & 1.034112e+00      \\ \hline
0.286032   & 1.589631e-05    & 1.023429e+00   & 1.023431e+00     & 1.627221e-03    & 1.034108e+00   & 1.034110e+00      \\ \hline
\end{tabular}}
\end{table}

\textbf{\textit{Order = 8}}
\begin{table}[h!]
\scalebox{0.9}{
\begin{tabular}{|c|c|c|c|c|c|c|}
\hline
\textbf{h} & \textbf{P comp} & \textbf{P err} & \textbf{Pmx err} & \textbf{U comp} & \textbf{U err} & \textbf{U mx err} \\ \hline
0.572063   & 3.409225e-05    & 1.023427e+00   & 1.023432e+00     & 2.404311e-03    & 1.034107e+00   & 1.034111e+00      \\ \hline
0.286032   & 1.037969e-05    & 1.023429e+00   & 1.023430e+00     & 1.349427e-03    & 1.034108e+00   & 1.034110e+00      \\ \hline
\end{tabular}}
\end{table}

\textbf{\textit{Order = 9}}
\begin{table}[h!]
\scalebox{0.9}{
\begin{tabular}{|c|c|c|c|c|c|c|}
\hline
\textbf{h} & \textbf{P comp} & \textbf{P err} & \textbf{Pmx err} & \textbf{U comp} & \textbf{U err} & \textbf{U mx err} \\ \hline
0.572063   & 2.328387e-05    & 1.023428e+00   & 1.023431e+00     & 2.033288e-03    & 1.034107e+00   & 1.034110e+00      \\ \hline
0.286032   & 7.124397e-06    & 1.023429e+00   & 1.023430e+00     & 1.141177e-03    & 1.034109e+00   & 1.034110e+00      \\ \hline
\end{tabular}}
\end{table}

\textbf{\textit{Order = 10}}
\begin{table}[h!]
\scalebox{0.9}{
\begin{tabular}{|c|c|c|c|c|c|c|}
\hline
\textbf{h} & \textbf{P comp} & \textbf{P err} & \textbf{Pmx err} & \textbf{U comp} & \textbf{U err} & \textbf{U mx err} \\ \hline
0.572063   & 1.664200e-05    & 1.023429e+00   & 1.023431e+00     & 1.746755e-03    & 1.034108e+00   & 1.034110e+00      \\ \hline
0.286032   & 5.085321e-06    & 1.023429e+00   & 1.023430e+00     & 9.803705e-04    & 1.034109e+00   & 1.034109e+00      \\ \hline
\end{tabular}}
\end{table}

\end{document}